\documentclass{article}
\usepackage{amsmath,amssymb,amsthm}

\newtheorem*{thmN}{Theorem}
\newtheorem{theorem}{Theorem}
\newtheorem{lemma}{Lemma}

\newtheorem{cor}{Corollary}
\theoremstyle{definition}
\newtheorem{definition}{Definition}
\newtheorem*{example}{Example}
\newtheorem*{rem}{Remark}
\newtheorem*{acknowledgments}{Acknowledgments}

\numberwithin{equation}{section}
\newcommand{\defn}{\ensuremath{\overset{\mathrm{def}}{=}}}
\newcommand{\dd}{\ensuremath{\mathrm{d}}}
\newcommand{\ddiv}{\ensuremath{\mathrm{Div}}}
\newcommand{\pr}{\ensuremath{\mathrm{pr}}}
\newcommand{\DD}{\ensuremath{\mathsf{D}}}

\newcommand{\sym}{\ensuremath{\mathrm{sym}}}
\newcommand{\dirac}{\ensuremath{\slash \!\!\! \partial}}

\begin{document}
\title{Conservation Laws and Non-Lie Symmetries}
\author{A.C.L Ashton\footnote{Department of Applied Mathematics and Theoretical Physics, University of Cambridge, Cambridge, CB3 0WA, UK. E-mail: a.c.l.ashton@damtp.cam.ac.uk}}
\maketitle

\begin{abstract}
\noindent
We introduce a method to construct conservation laws for a large class of linear partial differential equations. In contrast to the classical result of Noether, the conserved currents are generated by \emph{any} symmetry of the operator, including those of the non-Lie type. An explicit example is made of the Dirac equation were we use our construction to find a class of conservation laws associated with a 64 dimensional Lie algebra of discrete symmetries that includes CPT.
\end{abstract}

\section{Introduction}\label{intro}
The connection between symmetry and conservation laws has been inherent in all of mathematical physics since Emmy Noether published, in 1918, her hugely influential work linking the two \cite{noether1918iv}. Since then many have put forward approaches to study conservation laws, through a variety of different means \cite{anco2002dcm2,bluman2005cbs,finkel2001cee,fokas1996pib,gelfand1976fpo}. In each case, a conservation law is defined as follows.
\begin{definition}
Let $\Delta [u]=0$ be a system of equations depending on the independent variables $x=(x^1, \ldots , x^n)$, the dependent variables $u=(u^1, \ldots , u^m)$ and derivatives thereof. Then a \emph{conservation law} for $\Delta$ is defined by some $P=P[u]$ such that:
\begin{equation} \ddiv P\Big|_{\Delta =0} = 0 \label{divP} \end{equation}
where $[u]$ denotes the coordinates on the $N$-th jet of $u$, with $N$ arbitrary.
\end{definition}
Noether's theorem is applicable in the case where $\Delta[u]=0$ arises as the Euler-Lagrange equation to an associated variational problem. It is well known \cite{douglas1941sip,olver2000alg,anderson1992ivb,crampin1984gvh} that a PDE has a variational formulation if and only if it has self-adjoint Fr\'echet derivative. That is to say: if the system of equations $\Delta[u]=0$ is such that $\DD_\Delta = \DD_\Delta^*$ then the following result is applicable.
\begin{thmN}[Noether]
For a non-degenerate variational problem with $L[u]=\int_{\Omega}\mathfrak{L}\, \dd x$, the correspondence between nontrivial equivalence classes of variational symmetries of $L[u]$ and nontrivial equivalence classes of conservation laws is one-to-one.
\end{thmN}
Here the symmetries we speak of are \emph{generalised} symmetries which correspond to those with generators of the form:
\[ X = \sum_{i=1}^n \xi^i[u] \frac{\partial}{\partial x^i} + \sum_{i=1}^m \phi_i[u]\frac{\partial}{\partial u^i}. \]
Without loss of generality, one may confine attention to generalised vector fields in evolutionary form, where:
\[X_Q = \sum_{i=1}^m Q_i [u] \frac{\partial}{\partial u^i}. \]
The theorem of Noether then provides equivalence classes of conservation laws in terms of the characteristics of $X_Q$, i.e the $Q_i[u]$. It can be shown \cite{olver2000alg} that the $Q_i$ are related to the corresponding conservation laws via:
\begin{equation} Q\cdot \Delta = \ddiv P \label{char}\end{equation}
which holds for all $u$. This clearly gives rise to a conservation law since the LHS vanishes on solutions. Here we say $Q$ is the corresponding characteristic to the conservation law.

More generally, for systems that do not obey a variational formulation, we can consider \eqref{char} the definition of a conservation law. It can be shown \cite{olver2000alg} that the condition \eqref{divP} is equivalent to the existence of some $Q[u]$ such that $Q\cdot \Delta$ is a total divergence, as in \eqref{char}. Since the Euler operator annihilates divergences, this condition reads:
\begin{equation} E_i(Q\cdot \Delta)=0 \qquad 1 \leq i \leq m \label{chardef}\end{equation}
for each $(x,u)$. Here $E$ denotes the Euler operator:
\[ E_i = \sum_{\alpha}(-1)^{|\alpha|} D_\alpha \frac{\partial}{\partial u_\alpha^i}\]
where $\alpha$ is a multi-index and $u^i_\alpha \equiv \partial^\alpha u^i$. We regard \eqref{chardef} as the defining property for $Q$ to be a characteristic for a conservation law for $\Delta$. In this case, there is no interpretation of the $Q$ in terms of symmetries of the equations. Indeed, expanding the condition in \eqref{chardef} gives:
\[ \DD^*_\Delta (Q) + \DD_Q^* (\Delta) =0\]
for each $(x,u)$. In particular, on solutions we have what Anco and Bluman \cite{anco2002dcm1,anco2002dcm2} refer to as the \emph{adjoint-symmetry} condition:
\begin{equation} \DD_\Delta^* (Q)\Big |_{\Delta =0}=0. \label{adjsym} \end{equation}
This gives us a necessary, but not sufficient condition for $Q$ to be a characteristic for a conservation law. Now if $\Delta[u]=0$ corresponds to a variational problem, then $\DD_\Delta^* = \DD_\Delta$ and \eqref{adjsym} is automatically satisfied if $X_Q$ is the infinitesimal generator of a symmetry of $\Delta[u]=0$. This has an obvious connection to Noether's theorem: if $X_Q$ corresponds to variational symmetry of a variational problem, then it is also corresponds to a symmetry to the corresponding Euler-Lagrange equations (the converse is not true in general).

Throughout our discussion regarding symmetries, it has been implied that these symmetries correspond to a locally connected Lie group of transformations. If we were to talk of more general types of symmetry, then much of the previous material becomes void because of our use of infinitesimal generators. In addition, in our general definition of a conservation law, we insisted that it be a total divergence depending only on $[u]$ -- if we relax this, then it allows for a larger class of nontrivial total-divergences that vanish on the space of solutions. In this paper we address these points for a class of linear PDEs.

We consider a general set of symmetries, referred to as non-Lie symmetries \cite{nikitin1995nls,niederle1997nla,fushchich1971air}. This class includes generalized symmetries, non-local symmetries and discrete symmetries. An exposition into the construction of non-Lie symmetries can be found in \cite{nikitin1995nls} where several examples are given. Given that this class of symmetries is far larger than those considered in the classical work of Noether, there is potentially an even stronger correspondence between symmetry and conservation laws for PDEs. 

As an example, consider the scalar field $u(t,x)$, $x\in\Omega\subset\mathbf{R}^n$ whose evolution is prescribed via the heat equation, $u_t=\Delta u$. Then the following is a nontrivial conservation law for this system:
\[ D_t \left[ u^s u\right] + \ddiv \left[ u^s \nabla u - u\nabla u^s \right] =0 \] 
where $u^s(x,t) = u(x,s-t)$, $s$ fixed. Taking, for example, the solution to the related Cauchy problem with initial data $f\in C^\infty_c(\Omega)$ gives the family of conserved quantities
\[ E_s(t) = \frac{1}{[4\pi t(s-t)]^{n/2}} \int_{\mathbf{R}^{3n}}  e^{-\|x-y\|^2/4t - \|x-z\|^2/4(s-t)} f(y) f(z)\, \dd^n x\, \dd^ny\, \dd^n z \]
for $t\in (0,s)$. This result can not be achieved using standard techniques.

We formulate a method that provides insight into constructing new conservation laws generated from non-Lie symmetries, and that does not require a variational formulation. The method returns the classical results when the symmetries are of the Lie type, but provides \emph{new}, non-trivial conservation laws when algebras of discrete symmetries are used. We highlight this result by constructing new conservation laws for the Dirac equation, corresponding to a large Lie algebra of discrete symmetries. For example, we show that CPT symmetry in the Dirac equation generates the operator (after quantisation):
\[  \int \frac{\dd^3 p}{(2\pi)^3} \sum_s (-1)^s \left(  {a_{\mathbf{p}}^s}^\dagger a_{\mathbf{p}}^{s+1}+{b_{\mathbf{p}}^{s+1}}^\dagger b_{\mathbf{p}}^s \right), \]
which can be shown to commute with the Hamiltonian of the theory.

\section{Construction of the conservation law}
We work on a formal level so that all functions in question have sufficient smoothness and decay so that relevant integrals are well-defined. In particular, for evolution equations on $(0,T)\times \Omega \subset \mathbf{R}^{n+1}$ with conservation laws of the form:
\[ \frac{\partial \rho}{\partial t} + \ddiv J=0 \]
where $\rho=\rho[u]$ and $J=J[u]$, we assume the solutions are such that $\int_{\Omega}\rho\, \dd x$ converges and $J|_{\partial\Omega}=0$. Throughout we work with functions of $n+1$ independent variables so $C^\infty (\mathbf{R}^{n+1}, \mathbf{R}^m)$ etc. will be denoted simply by $C^\infty (\mathbf{R}^m)$. We concern ourselves with linear PDEs, so that in our previous notation $\Delta[u] \equiv \mathcal{L}[u]$ where $\mathcal{L}:C^{\infty}(\mathbf{R}^m)\rightarrow C^{\infty}(\mathbf{R}^l)$ is a linear matrix differential operator with smooth coefficients. We denote the formal adjoint of $\mathcal{L}$ by $\mathcal{L}^* : C^{\infty}(\mathbf{R}^l)\rightarrow C^{\infty}(\mathbf{R}^m)$, the coefficients of which are defined by $(\mathcal{L}^*)_{ij} = (\mathcal{L}_{ji})^*$, meaning that one first takes the transpose of the matrix differential operator $\mathcal{L}$, then takes the formal adjoint of each of differential operators. 

When we talk of \emph{symmetries} of a PDE, we refer to them in the sense of Fushchich and Nikitin \cite{fushchych1977nig}. The definition we use is as follows.
\begin{definition}\label{defsym}
We say the operator $\Gamma$ is a symmetry of the linear PDE $\Delta[u]\equiv \mathcal{L}[u]=0$ if there exists an operator $\mathcal{\alpha}_\Gamma$ such that:
\[ [\mathcal{L},\Gamma]=\alpha_\Gamma \mathcal{L} \]
where $[\cdot,\cdot]$ denotes the commutator by composition of operators so $\mathcal{L}\Gamma = \mathcal{L}\circ \Gamma$. We denote the set of all such symmetries by $\sym (\Delta)$.
\end{definition}
There is a rich source of literature on the study of such symmetries of linear PDOs \cite{fushchich1971air,fushchych1977nig,hydon2000cds,niederle1997nla} some of which we shall use in the sequel to construct conservation laws for the Dirac equation on Minkowski space. The analysis of such symmetries is not confined to those mentioned earlier, but also to those of the non-Lie type including discrete symmetries.
\begin{rem}
Definition \ref{defsym} can be viewed in terms of recursion operators. Recall that $\mathcal{R}$ is a recursion operator for the PDE $\Delta[u]=0$ if there exists an $\tilde{\mathcal{R}}$ such that:
\[ \DD_\Delta \mathcal{R} = \tilde{\mathcal{R}}\DD_\Delta. \]
In our case, $\Delta[u]$ is linear, so $\DD_\Delta$ coincides with $\mathcal{L}$. In terms of the notation in definition \ref{defsym} we have $\Gamma=\mathcal{R}$ and $\alpha_\Gamma = \tilde{\mathcal{R}}-\mathcal{R}$.

There is also a correspondence between these symmetries and the usual Lie-type as found in the standard treatment. Indeed, suppose the linear, $N$-th order PDE is invariant under the action of a one parameter group of transformations $\Gamma_\epsilon = \exp (\epsilon X)$, with corresponding infinitesimal generator:
\[ X = \sum_{i=1}^n \xi^i(x) \frac{\partial}{\partial x^i} + \sum_{i=1}^m \phi_i(x)\frac{\partial}{\partial u^i}.\]
Then by definition of $X$ we must have:
\[ \pr^N(X)\left(\Delta[u]\right)\Big|_{\Delta=0} = 0 \]
where $\pr^N(X)$ is the prolongation of $X$ onto $N$-th jet space. For this to hold, it must be the case that (see \cite{olver2000alg} p.82):
\[ \pr^N(X)\left(\Delta[u]\right) = \sum_{i=0}^M F_i(x) D_i \Delta[u], \]
for some $F_i(x)$ that are unimportant and some $M>0$. Now the LHS is equal to $\DD_\Delta (Xu)$, and since $\Delta$ is linear, this is just $\mathcal{L}[Xu]$. From this we see that the relevant commutator $[\mathcal{L},\pr^N(X)]$ takes the required form.
\end{rem}

From here on we deal entirely with symmetries of the form in definition \ref{defsym}. An immediate consequence of this definition is as follows.
\begin{lemma}\label{reclem}
For each $\Gamma \in \sym(\mathcal{L})$ and $u\in \ker (\mathcal{L})$ we have $\Gamma u\in \ker (\mathcal{L})$.
\end{lemma}
We saw in \S\ref{intro} that the existence of a conservation law to the equation $\Delta[u]=0$ was dependent on finding solutions to the so called adjoint-symmetry condition:
\begin{equation} \DD_\Delta^* (Q)\Big|_{\Delta=0}=0.\end{equation}
In the linear case with $\Delta[u]\equiv\mathcal{L}[u]$, this condition reads:
\begin{equation}
\mathcal{L}^*[Q]\Big|_{\Delta=0}=0 \label{symadj}
\end{equation}
The following lemma makes the condition also sufficient.
\begin{lemma}\label{symcons}
For a linear PDE $\Delta[u]\equiv \mathcal{L}[u]=0$, there exists a conservation law for each solution to \eqref{symadj}.
\end{lemma}
\begin{proof}
We present an argument similar to that presented in \cite{zharinov1986gsf}. Given $\mathcal{L}$, we associate with it the bilinear map $\Pi_{\mathcal{L}}: C^{\infty}(\mathbf{R}^m)\times C^{\infty}(\mathbf{R}^l) \rightarrow C^{\infty}(\mathbf{R})$, where $\Pi_\mathcal{L}$ is defined explicitly by:
\begin{equation} \Pi_\mathcal{L} (Q,P) = Q \cdot \mathcal{L}[P] - P\cdot \mathcal{L}^*[Q] \label{Pi}
\end{equation}
where $P,Q$ are $m$-tuples (respectively $l$-tuples) of smooth functions. Integration by parts shows this is a total divergence for each $P,Q$, so that $\Pi_{\mathcal{L}}(P,Q) \equiv \ddiv X$ for some $X$ depending bilinearly on $P$ and $Q$. Alternatively one can view this as a consequence of Green's formula \cite{zharinov1986gsf}. Setting $P=u$ we see that:
\[ \ddiv X|_{\Delta=0} = -u\cdot \mathcal{L}^*[Q]\Big|_{\Delta=0} \]
If $Q$ solves \eqref{symadj} then $\ddiv X=0$ on solutions to $\Delta [u]=0$, i.e a conservation law.
\end{proof}
\begin{cor}
If $\mathcal{L}$ is self-adjoint or skew-adjoint, then each $\Gamma \in \sym(\mathcal{L})$ generates a conservation law.
\end{cor}
The corollary follows by setting $Q=\Gamma u$ and appealing to the result in lemma \ref{reclem}.
\begin{example}
Consider the wave equation on $(0,T)\times\mathbf{R}^{n}$, $\mathcal{L} = \square$. Constructing the bilinear map as in the proof to lemma \ref{symcons} gives:
\[ \Pi_\mathcal{L} (Q,P) = D_t (Q P_t - PQ_t) + \ddiv \left(  P\nabla Q- Q\nabla P \right) \]
Setting $P=u$ and $Q=\Gamma u$ with $\Gamma \in \sym(\square)$ gives a conservation law. Integrating over $\mathbf{R}^n$ and the divergence theorem then gives the conserved quantities:
\[ \kappa_\Gamma = \int_{\mathbf{R}^n} \big [ (\Gamma u)u_t- uD_t (\Gamma u)\big ]\, \dd  x \]
for each symmetry $\Gamma$ of $\square$, assuming that all fields and their derivatives have sufficient decay at spatial $\infty$. Choosing $\Gamma = \partial_t$ gives rise to conservation of energy:
\begin{align*}
\kappa_{\partial_t} &= \int_{\mathbf{R}^n} \left( u_t^2 - uu_{tt}\right) \, \dd x \\
&= \int_{\mathbf{R}^n} \left( u_t^2 + |\nabla u|^2 \right) \, \dd x,
\end{align*}
which follows from integration by parts. This is in agreement with the classical result.
\end{example}
\begin{example}
We consider the coupled, linear KdV-KdV system \cite{bona2002bea,fokas2005bvp} on $(0,T)\times\mathbf{R}$:
\begin{align}
u_t + v_{xxx} + v_x &= 0, \label{lkdv1}\\
v_t + u_{xxx} + u_x &= 0.\label{lkdv2}
\end{align}
Writing the equations, $\Delta[u,v]=0$, in matrix form gives:
\[ \begin{bmatrix} D_t & D_x^3 + D_x \\ D_x^3 + D_{x} & D_t \end{bmatrix}\begin{bmatrix} u \\ v\end{bmatrix}=0 \]
We note that the system is skew-adjoint, and so the preceding results give a conservation law for each $\Gamma \in \sym(\Delta)$. Integrating this conservation law over $\mathbf{R}$ gives the conserved quantities:
\[ \kappa_\Gamma = \int_{\mathbf{R}} \big [ u (\Gamma u) + v(\Gamma v)\big ] \, \dd x \]
for each $\Gamma \in \sym (\Delta)$. Computing the symmetries by hand is straight forward, using the methods of Lie. The following vector fields correspond to the generators of Lie-point symmetries with linear coefficients:
\begin{align*} V_1 &= \partial_x, \quad V_2 = \partial_t, \quad V_3 = \partial_u, \quad V_4 = \partial_v, \\
\quad V_5 &= u\partial_v + v\partial_u, \quad V_6 = v\partial_u + u\partial_v, \quad
V_7 = t\partial_v -x\partial_u, \quad V_8 = x\partial_v - t\partial_u
\end{align*}
which generate an 8-dimensional Lie algebra. Labelling the corresponding conservation laws $\kappa_i$ we find $\kappa_1$, $\kappa_2$ are trivial and $\kappa_3,\kappa_4$ follow from the fact that \eqref{lkdv1}-\eqref{lkdv2} are already in conservation law form. Then $\kappa_5,\kappa_6$ corresponds to the conservation laws:
\[\kappa_5 = \int_{\mathbf{R}} (u^2 + v^2)\, \dd x,\quad \kappa_6 = \int_\mathbf{R} uv\, \dd x\]
Finally, the vector fields $V_7$, $V_8$ generate the conserved quantities:
\[ \kappa_7 = \int_{\mathbf{R}}\left( tu-xv\right)\, \dd x, \quad \kappa_8 = \int_{\mathbf{R}} \left( xu-tv\right)\, \dd x\]
Now we make use of a family of discrete symmetries defined by $\Gamma_s(u,v) = (-u^s, v^s)$ with $v^s(x,t) = v(x,s-t)$ etc, with $s$ fixed. Then for $t\in (0,s)$ we have the conserved quantity:
\[ \kappa_s = \int_\mathbf{R} \left[  vv^s -uu^s\right]\, \dd x. \]
To explicitly show this is indeed conserved, we compute the derivative.
\begin{align*}
\frac{\dd \kappa_s }{\dd t}  &= \int_\mathbf{R} \left[-(u_x + u_{xxx})v^s + (u^s_x + u^s_{xxx})v +(v_{xxx}+v_x)u^s - (v^s_{xxx}+v^s_x)u\right] \, \dd x\\
&= \int_{\mathbf{R}} D_x \left[ vu^s -uv^s + u^sv_{xx} + vu^s_{xx} - u^s_xv_x - v^su_{xx}- uv^s_{xx} + u_xv^s_x \right]\, \dd x \\
&=0
\end{align*}
which holds for each $t\in (0,s)$.
\end{example}
In this example we used an obvious discrete symmetry of the system. In general, to find the discrete symmetries of a system of PDEs we can apply the result of Hydon \cite{hydon2000cds}. In this method, one considers a general symmetry:
\[ \Gamma :(x,u)\rightarrow (\hat{x}(x,u),\hat{u}(x,u))\]
that may be continuous or discrete, then ``dresses'' this symmetry with a one-parameter group of symmetries $\Gamma_i(\epsilon)\equiv \exp(\epsilon X_i)$. The resulting symmetry:
\[ \hat{\Gamma}_i(\epsilon) = \Gamma \Gamma_i(\epsilon) \Gamma^{-1}=\exp(\epsilon \Gamma X_i \Gamma^{-1}) \]
is of the one parameter form, generated by $\hat{X}_i = \Gamma X_i \Gamma^{-1}$. The determining equations for $\hat{X}_i$ can then be related to the $X_i$, and the result is a system of nonlinear constraints that can be solved to give the most general form of symmetry. In particular, the discrete symmetries can be found by quotienting out the symmetries found via Lie's method. Several worked examples are given in \cite{hydon2000cds} which demonstrate the algorithmic process involved: the number of computations is not much more than that needed to compute the continuous symmetries.

We now turn our attention to the study of systems of PDEs with constant coefficients, and provide a general result that generates conserved quantities for each symmetry of the system. First we define the type of linear operators in question.
\begin{definition}
We denote the space of linear operators spanned by $M^\alpha D_\alpha$, $M_\alpha \in \mathrm{Mat}_m(\mathbf{C})$, by $L^c(\mathbf{C}^m)$. If $\mathcal{L}\in L^c(\mathbf{C}^m)$ then $\mathcal{L}:C^{\infty}(\mathbf{C}^m)\rightarrow C^\infty(\mathbf{C}^m)$.
\end{definition}
Such PDEs include both the wave equation and the linearised KdV-KdV system discussed earlier. Another example would be the linearised Navier-Stokes equations \cite{ashton2007fkf}:
\begin{align}
u_t +\nabla p - \nu \Delta u &=0 \label{ns1} \\
\ddiv (u) &=0 \label{ns2}
\end{align}
where $u=u(x,y,z,t)$ is the velocity field, $p=p(x,y,z,t)$ is pressure field and $\nu >0$. Then the corresponding element of $L^c(\mathbf{R}^4)$ is given by:
\[ \mathcal{L} = \begin{bmatrix} D_t - \nu\Delta & 0 & 0 & D_x \\
0 & D_t -\nu\Delta & 0 & D_y \\ 0 & 0 & D_t -\nu \Delta & D_z \\ D_x & D_y & D_x & 0 \end{bmatrix} \]
and equations \eqref{ns1}-\eqref{ns2} are given by $\mathcal{L}[\phi]=0$ with $\phi = (u^1, u^2, u^3, p)$.
\begin{definition}
We say $M^\alpha \in \mathrm{Mat}_m(\mathbf{C})$ are simultaneously semi-conjugate to $[M^\alpha]^\dagger$ if $\exists A_1, A_2 \in GL(m,\mathbf{C})$ independent of $\alpha$ such that $[M^\alpha]^\dagger = A_2 M^\alpha A_1^{-1}$ for each $\alpha$.
\end{definition}
\begin{lemma}\label{lem}
Suppose $\mathcal{L} = \sum_{\alpha} M^\alpha D_\alpha \in L^c(\mathbf{C}^m)$ and the $M^\alpha$ are simultaneously semi-conjugate to $[M^\alpha]^\dagger$. Then then there exists a pair of invertable linear operators, $\mathcal{P}_1, \mathcal{P}_2$, such that $\mathcal{L}^*=\mathcal{P}_2 \mathcal{L}\mathcal{P}_1^{-1}$.
\end{lemma}
\begin{proof}
Using $\mathcal{L}^* = \sum_\alpha (-1)^{|\alpha|}[M^\alpha]^\dagger D_\alpha$ we find:
\begin{equation}
\mathcal{L}^* = A_2 \left( \sum_\alpha (-1)^{|\alpha|} M^\alpha D_\alpha\right) A_1^{-1}. \label{lempr2}
\end{equation}
Now we introduce the parity operator $P:C^\infty(\mathbf{C}^m)\rightarrow C^{\infty}(\mathbf{C}^m)$ defined by $P\phi(x) = \phi(-x)$. We make the following observation:
\begin{align*}
(-1)^{|\alpha|} D_\alpha [\phi] &= P(D_\alpha) \phi(x) \\
&=P[D_\alpha (P\phi)] \\
&= (P\circ D_\alpha \circ P)[\phi].
\end{align*}
Using this identity in \eqref{lempr2} and the linearity of the associated operators we find:
\begin{equation}
\mathcal{L}^* = A_2 P \left( \sum_\alpha  M^\alpha D_\alpha\right)P A_1^{-1}.\label{lempr3}
\end{equation}
Then from the definition of $\mathcal{L}$ we see the required operators are $\mathcal{P}_2 = A_2 P$ and $\mathcal{P}_1 = PA_1^{-1} \equiv A_1^{-1}P$.
\end{proof}
Here we see the parity operator playing an important r\^ole, and it needs clarification as to whether this will cause any problems with our analysis. Terms such as $P[\phi]$ will be present in our conservation laws, and as such care must be taken when forming the corresponding conserved quantities.

In our opening remarks, we have made the assumption that in the case of evolution equations on $(0,T)\times \Omega \subset \mathbf{R}^{n+1}$ all fields are such that the relevant ``flux'', $J$, vanishes on $\partial\Omega$. In most cases this requires that the fields have sufficient decay at spatial $\infty$. If we say, for example, that $u(x,\cdot) \in W^{k,p}(\Omega)$, then it is clear that the discrete symmetries on the spatial independent variables will not cause any problems, assuming they are well-defined. For example if $\Omega = \mathbf{R}^n$, then if $u(x,\cdot) \in W^{k,p}(\Omega)$ then $u(-x,\cdot)\in W^{k,p}(\Omega)$ also.  Similarly if $\Omega$ is any isotropic domain. However, in the case of discrete symmetries of the temporal variables, they must take the form $\Gamma : u(\cdot, t) \mapsto u(\cdot, s-t)$, valid for $t\in (0,s)$ and $s \in (0,T)$. This means we need to modify our Parity operator so that when it acts on temporal variables, it does so in this prescribed form. As we shall see, however, in many instances the parity operator can be dispensed with altogether, by modifying the matrices $A_1, A_2$ in \eqref{lempr3} -- particularly in the case of Dirac operators (see lemma \ref{lem3} and corollary \ref{thmcor}).

In what follows we provide some lemmas that provide sufficient conditions for the matrices $M^\alpha$ to be simultaneously semi-conjugate to the $[M^\alpha]^\dagger$.
\begin{lemma}\label{lem2}
Suppose $\mathcal{L} = \sum_\alpha M^\alpha D_\alpha \in L^c(\mathbf{R}^m)$ such that $[M^\alpha, M^\beta]=0$ for each $\alpha, \beta$. Then $\mathcal{L}^* = \mathcal{P}_2 \mathcal{L} \mathcal{P}_1^{-1}$ for some linear, invertable $\mathcal{P}_1, \mathcal{P}_2$.
\end{lemma}
\begin{proof}
Since the $M^\alpha$ commute we may choose a basis $\{ \mathbf{e}_i\}_{i=1}^m$ such that the $M^\alpha$ simultaneously assume Jordan canonical form. This gives $M^\alpha = S_1 J^\alpha S^{-1}_1$ where $J^\alpha$ is in canonical form and $S_1\in GL(m,\mathbf{R})$. Now introduce a new basis $\{\tilde{\mathbf{e}}_i\}_{i=1}^m$ such that:
\[ \mathbf{e}_i \mapsto \tilde{\mathbf{e}}_i = \mathbf{e}_{m - i + 1}, \qquad 1\leq i \leq m \]
Under this change of basis, the $J^\alpha \mapsto [J^\alpha ]^t$. This then gives $M^\alpha = S_2 [J^\alpha]^t S_2^{-1}$ for some $S_2 \in GL(m,\mathbf{R})$. The result follows since $[J^\alpha]^t = S_3 [M^\alpha]^t S_3^{-1}$ for some $S_3 \in GL(m,\mathbf{R})$. Then $\mathcal{L}^* = \mathcal{P}_2 \mathcal{L} \mathcal{P}_1^{-1}$ for some linear, invertable $\mathcal{P}_1, \mathcal{P}_2$.
\end{proof}
\begin{lemma}\label{lem3}
Suppose $\mathcal{L} = \sum_\alpha M^\alpha D_\alpha \in L^c(\mathbf{R}^m)$ such that the $M_\alpha$ belong to a unitary representation of the Clifford algebra $Cl_{1,q}(\mathbf{R})$. Then $\mathcal{L}^* = \mathcal{P}_2 \mathcal{L} \mathcal{P}_1^{-1}$ for some linear, invertable $\mathcal{P}_1, \mathcal{P}_2$.
\end{lemma}
\begin{proof}
Recall that if $\{ \gamma^i\}$ are a unitary representation of $Cl_{1,q}(\mathbf{R})$ then they necessarily obey the anti-commutation relation:
\[ \gamma^i \gamma^j + \gamma^j \gamma^i= 2\eta^{ij} \mathbf{I}, \qquad \eta = \mathrm{diag} ( +1, \underbrace{-1,\ldots,-1}_q). \]
From this and the fact that $\gamma^i [\gamma^i]^\dagger = [\gamma^i]^\dagger \gamma^i = \mathbf{I}$ we find:
\[ [\gamma^i]^\dagger = \begin{cases} \gamma^i,& i=1 \\
                              -\gamma^{i},& 1< i\leq q+1
                             \end{cases}
\]
From this we deduce the result:
\[ [\gamma^i]^\dagger = \gamma^1 \gamma^i \gamma^1,\qquad 1\leq i \leq q+1. \]
So each of the $M^\alpha$ are simultaneously semi-conjugate to the $[M^\alpha]^\dagger$. The result now follows from lemma \ref{lem}.
\end{proof}
We now provide a concrete example to illustrate lemmas \ref{lem} and \ref{lem2}.
\begin{example}
Consider $\mathcal{L} = \sum_\alpha M^\alpha D_\alpha \in L^c(\mathbf{R}^2)$ given by:
\[ \mathcal{L} = \begin{bmatrix} D_t + D_x^3 & D_tD_x \\ 0 & D_t+D_x^3 \end{bmatrix} \equiv \begin{bmatrix} 1 & 0 \\ 0 & 1 \end{bmatrix} D_t + \begin{bmatrix} 1 & 0 \\ 0 & 1 \end{bmatrix}D_x^3 + \begin{bmatrix} 0 & 1 \\ 0 & 0 \end{bmatrix}D_t D_x. \]
Clearly each of the $M^\alpha$ commute, and since each matrix is already in Jordan canonical form we only need to compute the appropriate matrix $S_2\in GL(2,\mathbf{R})$ as in the proof to lemma \ref{lem2}. Using the result in \eqref{lempr3} gives:
\begin{align*} \mathcal{L}^* &= \begin{bmatrix} 0 & 1 \\ 1 & 0\end{bmatrix}P\left(  \begin{bmatrix} 1 & 0 \\ 0 & 1 \end{bmatrix} D_t + \begin{bmatrix} 1 & 0 \\ 0 & 1 \end{bmatrix}D_x^3 + \begin{bmatrix} 0 & 1 \\ 0 & 0 \end{bmatrix}D_t D_x\right)P\begin{bmatrix} 0 & 1 \\ 1 & 0\end{bmatrix} \\
&= -\begin{bmatrix} 1 & 0 \\ 0 & 1 \end{bmatrix} D_t - \begin{bmatrix} 1 & 0 \\ 0 & 1 \end{bmatrix}D_x^3 + \begin{bmatrix} 0 & 0 \\ 1 & 0 \end{bmatrix}D_t D_x \\
&= \begin{bmatrix} -D_t - D_x^3 & 0 \\ D_t D_x & -D_t-D_x^3 \end{bmatrix}
\end{align*} 
Which agrees with the expected result.
\end{example}
We note that there exists an alternate argument for the result in Lemma \ref{lem} which applies to \emph{all} elements in $L^c(\mathbf{R}^m)$. Indeed, if we consider $\mathcal{L}$ as a matrix over the ring of pseudo-differential operators \cite{dickey1990sea,shubin1987poa} of the form:
\begin{equation} X_\Psi =  \sum_{-\infty}^{\mathrm{finite}} c^\alpha(\phi) D_\alpha \label{pseudo} \end{equation}
we have considerable freedom. In our case, the $c^\alpha$ are constants and as such the associated pseudo-differential operators form a commutative ring. What is more, inverses exist -- so on a formal level we may treat our $\mathcal{L}$ as matrices over this field of differential operators. Then recalling from linear algebra that each matrix over a field $F$ is similar to its transpose, we may use a similar argument to that found in lemma \ref{lem} to prove the result for arbitrary $\mathcal{L}\in L^c(\mathbf{R})$. However, as we shall see, we will be interested on the action of the relevant operators $\mathcal{P}_i$ on test functions $\phi \in C^{\infty}(\mathbf{R}^m)$, so the inclusion of operators of the form in \eqref{pseudo} causes considerable technical issues of how one should define $X_\Psi [\phi]$. For this reason we treat the indicated subset of $L^c(\mathbf{R}^m)$.
\begin{theorem}\label{consvlaw}
Let $\Delta[u]=0$ be a linear system of PDEs associated with the operator $\mathcal{L} =\sum_\alpha M^\alpha D_\alpha\in L^c(\mathbf{R}^m)$, where the $M^\alpha$ are simultaneously semi-conjugate to the $[M^\alpha]^\dagger$. Then for each $\Gamma \in \sym (\Delta)$ there exists a conservation law.
\end{theorem}
\begin{cor}\label{thmcor}
Suppose the $\mathcal{L}$ associated with $\Delta[u]=0$ is a polynomial in the Dirac operator $\dirac  \equiv \sum_i \gamma^iD_i$. Then for each $\Gamma \in \sym (\Delta)$ there exists a conservation law.
\end{cor}
\begin{proof}
We see from lemma \ref{reclem} that for each $\Gamma \in \sym (\Delta)$ we have $\Gamma u \in \ker (\mathcal{L})$ for solutions to $\Delta [u]=0$. Also by lemma \ref{lem} it follows that if $u\in \ker (\mathcal{L})$ then $\mathcal{P}^{-1}_1u \in \ker( \mathcal{L^*})$. Indeed, using the result from lemma \ref{lem} we have:
\[ \mathcal{L}^* \mathcal{P}_1^{-1}u = \mathcal{P}_2 \mathcal{L} \mathcal{P}_1 \mathcal{P}^{-1}_1 u = \mathcal{P}_2 \mathcal{L}u =0 \]
since $u\in\ker (\mathcal{L})$. Combining these results gives:
\[ u \in \ker (\mathcal{L}) \quad \Rightarrow \quad \mathcal{P}^{-1} \Gamma u \in \ker (\mathcal{L}^*) \]
for each $\Gamma \in \sym (\Delta)$. This can be restated as:
\[ \mathcal{L}^* \left( \mathcal{P}^{-1} \Gamma u\right) \Big|_{\Delta = 0} = 0 \]
for each $\Gamma \in \sym (\Delta)$. The result now follows after an application of lemma \ref{symcons}. The corollary follows from this result and an application of lemmas \ref{lem} and \ref{lem3}.
\end{proof}
We give an extended example in the next section which utilizes the results in theorem \ref{consvlaw} and corollary \ref{thmcor}. 

\section{The Dirac Equation}
Dirac's celebrated equation governing spin $\tfrac{1}{2}$ fermions takes the covariant form:
\begin{equation} (i\dirac -m)\psi = 0, \qquad \dirac \equiv \gamma^0 D_t - \gamma^i D_i \label{dirac_eqn} \end{equation}
where here and throughout this section we employ summation convention (roman letters running over 1,2,3 and Latin over 0,1,2,3). The $\{\gamma^\mu\}$ are members of the Dirac algebra $Cl_{1,3}(\mathbf{C})=Cl_{1,3}(\mathbf{R})\otimes \mathbf{C}$ which obey:
\begin{equation} \{ \gamma^\mu, \gamma^\nu\} \defn \gamma^\mu \gamma^\nu + \gamma^\nu \gamma^\mu = 2\eta^{\mu\nu}\mathbf{I} \label{cliffrel} \end{equation}
where $\eta$ is the metric on Minkowski space, with signature $(+,-,-,-)$. Here $\psi$ is a Dirac spinor associated with the $(\frac{1}{2},0)\oplus(0,\frac{1}{2})$ representation of the Lorentz group. We employ the Dirac representation of $Cl_{1,3}(\mathbf{C})$, so that the $\{\gamma^\mu\}$ are given by the following $4\times 4$ matrices:
\begin{equation}
\gamma^0 = \begin{pmatrix}  \mathbf{I}_2 & 0 \\ 0& -\mathbf{I}_2  \end{pmatrix}, \quad \gamma^i = \begin{pmatrix} 0 & \sigma^i \\ -\sigma^i & 0 \end{pmatrix}.
\end{equation}
where the $\sigma^i$ are the Pauli spin matrices. We use $\eta^{\mu\nu}$ to raise and lower indices so that $\gamma_\mu = \eta_{\mu \nu} \gamma^\nu$ where $\eta_{\mu \nu} \eta^{\nu \rho} = \delta_\mu^\rho$. The field equation (\ref{dirac_eqn}) comes from varying the following action:
\[ \mathcal{S}[\psi,\bar{\psi}]=\int \bar{\psi} (i \dirac - m)\psi\, \dd x. \] 
Here $\bar{\psi} \equiv \psi^\dagger \gamma^0$ is referred to as the Dirac conjugate. This variational formulation admits a $U(1)$ variational symmetry\footnote{Or gauge symmetry in the physics literature.}, from which we deduce the conserved current $j^\mu = \bar{\psi}\gamma^\mu \psi$ via Noether's theorem. We will now derive this result without referring to variational principles. We identify the appropriate $\mathcal{L} \in L^c(\mathbf{R}^4)$ as:
\begin{equation}
\mathcal{L} = \begin{bmatrix} (iD_t - m)\mathbf{I}_2 & -i\sigma^i D_i \\ i\sigma^iD_i & -(iD_t + m)\mathbf{I}_2 \end{bmatrix}.
\end{equation}
Now since this is a (linear) polynomial of Dirac operators, corollary \ref{thmcor} applies and we may construct a conservation law from each symmetry of \eqref{dirac_eqn}. Noting that the fields are complex, we construct $\Pi_\mathcal{L}$ using the adjoint associated with the $L^2$ inner product $(f,g) = \int f g^*\, \dd x$. A routine calculation gives:
\begin{align*}
 \Pi_\mathcal{L}(\psi, \tilde{\psi}) &= \tilde{\psi}_i^* (i\dirac -m )_{ij} \psi_j + \psi_j (i\dirac + m)_{ij} \tilde{\psi}_i^* \\
&= \partial_\mu (i\tilde{\psi}^\dagger \gamma^\mu \psi ).
\end{align*}
Since $\mathcal{L}$ is a polynomial in $\dirac$, the result in lemma \ref{lem3} gives:
\begin{equation}\mathcal{L}^* =  \gamma^0 (i\dirac - m)\gamma^0 \end{equation}
Then in the notation of lemma \ref{lem}, we have that $\mathcal{P}_1 = \mathcal{P}_2= \gamma^0$. Now using the result of theorem \ref{consvlaw} with $\Gamma = \mathrm{id} \in \sym(\Delta)$ gives the conservation law:
\begin{align*} \Pi_\mathcal{L}(\psi, \psi) &= \partial_\mu (i(\gamma^0 \psi)^\dagger \gamma^\mu \psi) \\
 &= \partial_\mu (i\psi^\dagger \gamma^0 \gamma^\mu \psi) \\
 &= \partial_\mu (i \bar{\psi} \gamma^\mu \psi),
\end{align*}
agreeing with the classical result corresponding to conservation of charge. However, we may now use the result from the previous section to create more conservation laws. If $\Gamma \in \sym (\Delta)$ then by theorem \ref{consvlaw}:
\begin{equation}\Pi_\mathcal{L} (\psi, \gamma^0\Gamma \psi)=0 \label{cur}\end{equation}
is a conservation law associated with the Dirac equation. Each of these generate a conserved quantity:
\begin{equation}
\kappa_\Gamma = \int (\Gamma \psi)^\dagger \psi\, \dd x \label{dens}
\end{equation} 
which follows from \eqref{cur} and $(\gamma^0)^2=\mathbf{I}$. This holds for \emph{all} symmetries of \eqref{dirac_eqn}, not only those corresponding to variational symmetries in the Lagrangian formulation. For completeness, we observe that the result in \eqref{dens} agrees with the classical results when $\Gamma \in \sym (\Delta)$ is a symmetry of the Lie type.

It is well known that the maximal algebra of Lie type symmetries of the Dirac operator is the Poincar\'{e} algebra, which has basis:
\[ P_\mu = -iD_\mu, \quad J_{\mu\nu} = -i (x_\mu D_\nu - x_\nu D_\mu) + \tfrac{i}{4} [\gamma_\mu, \gamma_\nu]. \]
We focus on the symmetries that generate rotations, i.e:
\[ L_i := \tfrac{1}{2}\epsilon_{ijk}J_{jk}. \]
It is a simple exercise to show that this gives rise, by means of Noether's theorem, to conservation of angular momentum (orbital and spin) $\mathbf{J}$ where:
\begin{equation}
\mathbf{J} =  \int  \psi^\dagger \left(\mathbf{x}\wedge\mathbf{p}  + \tfrac{1}{2} \mathbf{\Sigma}\right)\psi\, \dd x  \label{ang_mom}
\end{equation}
where $\mathbf{p}=-i\nabla$ and $\mathbf{\Sigma} = \mathrm{diag} (\sigma_i, \sigma_i)$. We use the same subset of the Poincar\'{e} algebra that generates rotations, so the associated symmetry operator is:
\begin{align*} \Gamma_i &= \tfrac{1}{2}\epsilon_{ijk}\left(x_jp_k - x_kp_j + \tfrac{i}{4}[\gamma_i, \gamma_j]\right)\\
&= \epsilon_{ijk} x_j p_k + \tfrac{1}{4}\epsilon_{ijk}\epsilon_{jkl}\Sigma_l \\
&= \epsilon_{ijk} x_j p_k + \tfrac{1}{2}\Sigma_i
\end{align*}
We then see this leads to the same conserved quantity as derived in \eqref{ang_mom}. Clearly each of the conserved quantities that follow from Noether's theorem can also be generated by this method, but the important point to outline is that we may construct other conserved quantities, generated by non-Lie symmetries, including discrete.

In \cite{niederle1997nla} the algebraic aspects of the discrete symmetries of the Dirac equation were studied in detail. We follow this work, and study the algebra associated with the discrete parity, charge and time symmetries present in the Dirac equation. The basis for this algebra is:
\[ \Gamma_\mu = \gamma_4\gamma_\mu \hat{\theta}_\mu, \quad \Gamma_4 = i\gamma_4 \hat{\theta}, \quad \Gamma_5 = i\gamma_2 \hat{c},\quad \Gamma_6 = i\Gamma_5 \]
where $\hat{\theta}_2 \psi(x) = \psi (x_0, x_1, -x_2, x_3)$ etc, $\hat{\theta}\psi (x) = \psi(-x)$, $\hat{c}\psi = \psi^*$ and $\gamma_4$ is the fourth\footnote{Or often referred to as the fifth.} gamma matrix:
\[ \gamma_4 = i\gamma_0\gamma_1\gamma_2\gamma_3. \]
These basis elements actually correspond to a 7-dimensional Clifford algebra:
\[ \{\Gamma_a, \Gamma_b \} = g_{ab} \mathbf{I},\qquad g = \mathrm{diag} (+1,-1,-1,-1,+1,+1,+1).\]
It can be shown that the enveloping algebra of this Clifford algebra is isomorphic to $\mathfrak{gl}(8,\mathbf{R})$ which means there are 64 linearly independent combinations of these elements. The algebraic aspects of these symmetries has a rich structure, but we shall not pursue it here so refer the reader to \cite{fushchich1971air,fushchych1977nig,niederle1997nla}. We now form a conserved quantity from this an element of this algebra: we choose $\Gamma = \Gamma_0\in \sym (\Delta)$ for simplicity. It follows from \eqref{dens} that the following quantity is conserved:
\[ \kappa_0 =  \int \bar{\psi}(-t,\mathbf{x})\gamma_4 \psi (t,\mathbf{x})\, \dd x. \]
Noting that $\bar{\psi}(-t,\mathbf{x})$ satisfies $i\partial_0 \bar{\psi}\gamma^0 = -i\partial_i\bar{\psi}\gamma^i + m\bar{\psi}$, we find:
\begin{align*}
\frac{\dd \kappa_0}{\dd t} &= \int  \left[ i \partial_i \bar{\psi}(-t,\mathbf{x})\gamma^i \gamma^4 \psi (t,\mathbf{x}) - i\bar{\psi}(-t,\mathbf{x})\gamma^4\gamma^i \partial_i \psi (t,\mathbf{x}) \right ]\,\dd x \\
&= \int \left[ i \partial_i \bar{\psi}(-t,\mathbf{x})\gamma^i \gamma^4 \psi (t,\mathbf{x}) + i\bar{\psi}(-t,\mathbf{x})\gamma^i\gamma^4 \partial_i \psi (t,\mathbf{x}) \right ]\,\dd x \\
&= \int  \partial_i \left[ \bar{\psi}(-t,\mathbf{x})\gamma^i\gamma^4 \psi(t,\mathbf{x})\right]\,\dd x \\
&=0.
\end{align*}
Where we used the fact that $\gamma_4$ anti-commutes with the $\gamma_\mu$ and the assumption that the fields decay at spatial $\infty$. We can express this quantity in terms of the usual creation and annihilation operators -- a short calculation reveals:
\begin{eqnarray}
\kappa_0 =\int \frac{\dd^3 p}{(2\pi)^3} \, \sum_s \left( {a_{-\mathbf{p}}^s}^\dagger b_{\mathbf{p}}^s + {b_{-\mathbf{p}}^s}^\dagger a_{\mathbf{p}}^s \right) \label{kappa_0}
\end{eqnarray}
where we have employed the usual quantisation procedure of promoting the fields to operators on a Fock space, and the $a^s_\mathbf{p}$, $b^s_{\mathbf{b}}$ etc. obey the anti-commutation relations:
\[ \{a^s_\mathbf{p}, {a^r_\mathbf{q}}^\dagger\} = \{b^s_\mathbf{p}, {b^r_\mathbf{q}}^\dagger\}=(2\pi)^3 \delta^{rs} \delta^3(\mathbf{p}-\mathbf{q}) \]
in the distributional sense, and all other anti-commutators being zero. We see $\kappa_0$ annihilates the vacuum, and using the usual anti-commutation relations we find:
\[ [\kappa_0,{a^s_{\mathbf{p}}}^\dagger] =  {b^s_{-\mathbf{p}}}^\dagger, \quad [\kappa_0,{b^s_{\mathbf{p}}}^\dagger] =  {a^s_{-\mathbf{p}}}^\dagger \]
So $\kappa_0$ swaps an anti-particle for a particle, and reverses the direction of momentum. A routine calculation shows that $\kappa_0$ commutes with the Hamiltonian $H$:
\[ [ H, \kappa_0] = 0, \quad H= \int \frac{\dd^3 p}{(2\pi)^3} \, E_\mathbf{p} \sum_s \left( {a_{\mathbf{p}}^s}^\dagger a_{\mathbf{p}}^s + {b_{\mathbf{p}}^s}^\dagger b_{\mathbf{p}}^s \right) \]
so it is indeed conserved. It is interesting to note the conserved quantity associated with the \emph{CPT} symmetry that is inherent in all field theories governed by a Hermitian Hamiltonian. In this case it is given by $\Gamma_4\Gamma_5$, and the corresponding conserved quantity is:
\begin{equation}
\kappa_{45} = \int \dd^3 x\, \bar{\psi}(-t,-\mathbf{x})^* \gamma^2\gamma^0\gamma^4 \psi (t,\mathbf{x}), \label{kappa_45}
\end{equation}
which is easily shown to be constant in time. Expressing this in terms of the creation and annihilation operators for the particles and anti-particles, we find that the operator in (\ref{kappa_45}) is equivalent to:
\begin{equation}
\int \frac{\dd^3 p}{(2\pi)^3} \sum_s (-1)^s \left(  {a_{\mathbf{p}}^s}^\dagger a_{\mathbf{p}}^{s+1}+{b_{\mathbf{p}}^{s+1}}^\dagger b_{\mathbf{p}}^s \right), \label{kappa_45'}
\end{equation}
where the spin sum is over $\mathbf{Z}/2\mathbf{Z}$. The derivation of this result can be found in the appendix. A routine calculation confirms that this operator commutes with the Hamiltonian.

\section{Conclusions}
In this paper we have seen that the number of conservation laws for a system of linear PDEs can be greatly increased if one considers the associated densities to depend not only on $[u]$, but also \emph{discrete} transformations of thereof. Allowing for this extra freedom, we have produced a general result applicable to a large class of linear PDEs, including those associated with polynomials of Dirac operators. The results of Hydon \cite{hydon2000cds} and Fushchich \cite{fushchich1971air} allows for an algorithmic construction of the non-Lie symmetries, upon finding which theorem \ref{consvlaw} provides us with a means to construct corresponding conservation laws.

\begin{acknowledgments}
The author would like to thank Thanasis Fokas for his lessons in perseverance and George Weatherill for many useful discussions. In addition, the author thanks the anonymous referee who made him aware of the work in \cite{zharinov1986gsf}. This work was supported EPSRC.
\end{acknowledgments}

\section{Appendix}
\subsection*{$CPT$ Symmetry in the Dirac Equation}
It was shown in \S 3 that the $CPT$ symmetry inherent in the Dirac equation gives rise to the conserved quantity:
\begin{equation}
\int \frac{\dd^3 p}{(2\pi)^3} \sum_s (-1)^s \left( {a_{\mathbf{p}}^s}^\dagger a_{\mathbf{p}}^{s+1} + {b_{\mathbf{p}}^{s+1}}^\dagger b_{\mathbf{p}}^s \right), \label{appen1}
\end{equation}
where the spin sum is over $\mathbf{Z}/2\mathbf{Z}$. In this section we outline the details of the calculation that produces this result, after applying the standard quantisation procedure to the term:
\begin{equation} \int \dd^3 x\, \bar{\psi}(-t,-\mathbf{x})^* \gamma^2\gamma^0\gamma^4 \psi (t,\mathbf{x}), \label{appen2} \end{equation}
The calculations in this section would be somewhat easier if we chose to work with Marjona spinors from the offset, but we shall stick with our initial choice and suffer the consequences. Firstly we make an observation:
\begin{align*}
 & 0 = (i\dirac - m)\psi (x), \\
\Rightarrow\quad &0 = (-i \left[\gamma ^*\right]^\mu \partial_\mu - m) \psi (x)^* , \\
\Rightarrow\quad &0 = (+i\gamma^0\partial_t - i\gamma^1 \partial_x + i\gamma^2\partial_y - i\gamma^3\partial_z - m) \psi (x)^* , \\
\Rightarrow\quad &0 = (-i\gamma^0\partial_t + i\gamma^1 \partial_x + i\gamma^2\partial_y + i\gamma^3\partial_z - m) \psi (-t,-x,y,-z)^* , \\
\Rightarrow\quad &0 = (i\dirac - m)\psi (-t,-x,y,-z)^*
\end{align*}
And so we set:
\begin{equation}
\bar{\psi}(-t,-\mathbf{x})^* = \bar{\psi}(t,x,-y,z).
\end{equation}
We shall use this in \eqref{appen2}, along with the mode expansions for $\psi$, $\bar{\psi}$. In the Sch\"{o}dinger picture the relevant mode expansions are given by:
\begin{align*}
\psi (x,y,z) &= \int \frac{\dd^3 p}{(2\pi)^3} \frac{1}{\sqrt{2E_{\mathbf{p}}}} \sum_s \left( a^s_{\mathbf{p}} u_s(\mathbf{p}) + {b^s_{-\mathbf{p}}}^\dagger v_s(-\mathbf{p})\right) e^{-i\mathbf{p}\cdot\mathbf{x}} \\
\bar{\psi} (x,-y,z) &= \int \frac{\dd^3 p}{(2\pi)^3} \frac{1}{\sqrt{2E_{\mathbf{p}}}} \sum_s \left( {a^s_{\mathbf{p}}}^\dagger \bar{u}_s(\mathbf{p}) +b^s_{-\mathbf{p}} \bar{v}_s(-\mathbf{p})\right) e^{i\mathbf{p}'\cdot\mathbf{x}}
\end{align*}
where $\mathbf{p}' = (p^1, -p^2, p^3)$. Using these expressions in \eqref{appen2} and doing the $x$ integral and one of the $\mathbf p$ integrals, we find:
\begin{equation}
\int \frac{\dd^3 p}{(2\pi)^3} \frac{1}{2E_{\mathbf{p}}} \sum_{r,s} \left( {a^r_{\mathbf{p}}}^\dagger \bar{u}_r(\mathbf{p}) +b^r_{-\mathbf{p}} \bar{v}_r(-\mathbf{p})\right) \gamma^2 \gamma^0\gamma^4 \left( a^s_{\mathbf{p}'} u_s(\mathbf{p}') + {b^s_{-\mathbf{p}'}}^\dagger v_s(-\mathbf{p}')\right). \label{big}
\end{equation}
Using the anti-commutation relation for the gamma matrices in \eqref{cliffrel}, we have $\gamma^2 \gamma^0 \gamma^4 = \gamma^0 \gamma^4 \gamma^2$. Inserting the $\gamma^2$ matrix into the right hand side bracket in \eqref{big} leaves us to compute $\gamma^2 u_s (\mathbf{p}')$ and $\gamma^2 v_s (-\mathbf{p}')$. We consider $\gamma^2 u_s (\mathbf{p}')$ only, since the computation for $\gamma^2 v_s (-\mathbf{p}')$ is entirely analogous. In our representation, $\gamma^2 u_s (\mathbf{p}')$ is given by:
\[ \sqrt{E_{\mathbf{p}'} +m} \begin{pmatrix} 0 & \sigma^2 \\ -\sigma^2 & 0 \end{pmatrix} \begin{pmatrix} \chi_s \\ \frac{(\mathbf{\sigma}\cdot\mathbf{p}')\chi_s}{E_{\mathbf{p}'}+m} \end{pmatrix}, \]
where $\chi_1 =\left[\begin{smallmatrix} 1 \\ 0\end{smallmatrix}\right]$ and $\chi_2 = \left [\begin{smallmatrix} 0 \\ 1\end{smallmatrix}\right]$ are 2-spinors. Using the fact that $E_{\mathbf{p}'} = E_{\mathbf{p}}$ and $\sigma^2 (\mathbf{\sigma}\cdot\mathbf{p}') = - (\mathbf{\sigma}\cdot \mathbf{p})\sigma^2$ we find:
\[ \gamma^2 u_s (\mathbf{p}') = i(-1)^s v_{s+1}(\mathbf{p}), \qquad \mathrm{where} \,\,\, v_s(\mathbf{p}) = \sqrt{E_{\mathbf{p}'} +m}\begin{pmatrix} \frac{(\mathbf{\sigma}\cdot\mathbf{p})\chi_s}{E_{\mathbf{p}}+m}  \\ \chi_s \end{pmatrix} \]
where $s\in \mathbf{Z}/2\mathbf{Z}$. A similar computation for $\gamma^2 v_s (-\mathbf{p}')$ reveals:
\[ \gamma^2 v_s (-\mathbf{p}') = i(-1)^{s+1} u_{s+1}(\mathbf{p}). \]
The expression in \eqref{big} now takes the form:
\begin{multline}\int \frac{\dd^3 p}{(2\pi)^3} \frac{i}{2E_{\mathbf{p}}} \sum_{r,s} \left( {a^r_{\mathbf{p}}}^\dagger \bar{u}_r(\mathbf{p}) +b^r_{-\mathbf{p}} \bar{v}_r(-\mathbf{p})\right) \\
\times \gamma^0\gamma^4 \left( a^s_{\mathbf{p}'} (-1)^s v_{s+1}(\mathbf{p}) + {b^s_{-\mathbf{p}'}}^\dagger (-1)^{s+1} u_{s+1}(\mathbf{p})\right).\label{bigger} \end{multline}
Now we must consider terms of the form $\bar{u}_r (\mathbf{p}) \gamma^0\gamma^4 v_{s+1}(\mathbf{p})$ and similar. We perform the calculation for $\bar{u}_r (\mathbf{p}) \gamma^0\gamma^4 v_{s+1}(\mathbf{p})$ since the computation is similar for the remaining three terms. In our representation the computation reads:
\begin{align*}
\bar{u}_r (\mathbf{p})\left( \begin{smallmatrix} 0 & \mathbf{I}_2 \\ -\mathbf{I}_2 & 0 \end{smallmatrix}\right) v_{s+1}(\mathbf{p}) &= (E_{\mathbf{p}}+m) \left( \chi_r^\dagger, -\chi_r^\dagger \left [\tfrac{(\sigma\cdot\mathbf{p})}{E_\mathbf{p}+m}\right] \right) \begin{pmatrix} \chi_{s+1} \\ - \left [\tfrac{(\sigma\cdot\mathbf{p})}{E_\mathbf{p}+m}\right]\chi_{s+1} \end{pmatrix} \\
&= (E_\mathbf{p}+m )\left [ \chi_r^\dagger \chi_{s+1} + \chi_r^\dagger \left [\tfrac{(\sigma\cdot\mathbf{p})}{E_\mathbf{p}+m}\right]^2 \chi_{s+1} \right ] \\
&= \left[ E_\mathbf{p}+m + E_\mathbf{p} -m \right] \delta^{r,s+1} \\
&= 2E_\mathbf{p} \delta^{r,s+1}
\end{align*}
in which we used $(\sigma \cdot \mathbf{p})^2 = |\mathbf{p}|^2 = E_\mathbf{p}^2 - m^2$. A similar calculation reveals:
\[ \bar{v}_r(-\mathbf{p}) \gamma^0\gamma^4 u_{s+1}(\mathbf{p}) = 2E_\mathbf{p} \delta^{r,s+1},\]
with the remaining two terms in \eqref{bigger} evaluating to zero. Using these results we find that \eqref{bigger} reduces to:
\begin{equation}
\int \frac{\dd^3 p}{(2\pi)^3} \sum_s i(-1)^s \left( { a_{\mathbf{p}}^s}^\dagger a_{\mathbf{p}}^{s+1} -  b_{-\mathbf{p}}^s {b_{-\mathbf{p}}^{s+1}}^\dagger \right).
\end{equation}
Now discarding the unimportant constant $i$, changing $\mathbf{p}\mapsto -\mathbf{p}$ in the second term and using the anti-commutivity of the $b$'s, we recover the expression in \eqref{appen1}.

\nocite{ashton2007fkf}
\nocite{finkel2001cee}
\nocite{minzoni2000scq}
\small
\bibliographystyle{plain}
\bibliography{../../../Bibliography/ants_bib}

\label{lastpage}
\end{document}